\documentclass[11pt,reqno,a4paper]{amsart}

\usepackage{amsfonts, amsmath, amsthm, amssymb}
\usepackage[english]{babel}
\usepackage[mathcal]{eucal}

\allowdisplaybreaks

\title[]{An approach theoretic version of Anscombe's theorem with an application in biostatistics}
\author{Ben Berckmoes}
\dedicatory{Dedicated to the 70th birthday of Bob Lowen, my academic father, an outstanding mathematician, an extraordinary teacher, and a beautiful human being.}
\thanks{Ben Berckmoes is post doctoral fellow at the Fund for Scientific Research of Flanders (FWO)}

\date{}

\DeclareMathOperator*{\myinf}{in\vphantom{p}f}

\begin{document}

\maketitle

\newtheorem{pro}{Proposition}[section]
\newtheorem{lem}[pro]{Lemma}
\newtheorem{thm}[pro]{Theorem}
\newtheorem{de}[pro]{Definition}
\newtheorem{co}[pro]{Comment}
\newtheorem{no}[pro]{Notation}
\newtheorem{vb}[pro]{Example}
\newtheorem{vbn}[pro]{Examples}
\newtheorem{gev}[pro]{Corollary}
\newtheorem{vrg}[pro]{Question}
\newtheorem{rem}[pro]{Remark}
\newtheorem{lemA}{Lemma}

\begin{abstract}
We establish an approach theoretic version of Anscombe's theorem, which we apply to justify the use of confidence intervals based on the sample mean after a group sequential trial.
\end{abstract}

\section{Introduction}

In Anscombe theory, the stability of weak convergence of random variables under index randomization is investigated. Let us explain this idea more precisely.

Let $X$ be a separable metric space with metric $d$, and consider $X$-valued random variables $\xi$ and $(\xi_n)_n$. We let $\stackrel{w}{\rightarrow}$ stand for weak convergence, i.e.\ $\xi_n \stackrel{w}{\rightarrow} \xi$ means that, for all $f : X \rightarrow \mathbb{R}$ bounded and continuous, 
$$\mathbb{E}[f(\xi_n)]\rightarrow \mathbb{E}[f(\xi)].$$ 
Furthermore, consider a sequence $(N_n)_n$ of $\mathbb{N}$-valued random variables. Anscombe theory deals with the following natural question. Under which conditions does the implication 
$$\xi_n \stackrel{w}{\rightarrow} \xi \Rightarrow \xi_{N_n} \stackrel{w}{\rightarrow} \xi$$ 
hold?

Let $\stackrel{\mathbb{P}}{\rightarrow}$ stand for convergence in probability, i.e.\ $\xi_n \stackrel{\mathbb{P}}{\rightarrow} \xi$ means that, for each $\epsilon > 0$,
$$\mathbb{P}[d(\xi,\xi_n) \geq \epsilon] \rightarrow 0.$$ 
Furthermore, we say that a sequence $(\xi_n)_n$ of random variables satisfies Anscombe's condition iff for each $\epsilon > 0$ there exist $\delta > 0$ and $n_0$ such that 
$$\mathbb{P}\left[\max_{m = \lceil (1 - \delta) n \rceil}^{\lfloor (1 + \delta) n \rfloor} d(\xi_n,\xi_m) \geq \epsilon\right] < \epsilon$$ 
for all $n \geq n_0$, $\lfloor \cdot \rfloor$ being the floor function, and $\lceil\cdot\rceil$ the ceiling function. This condition guarantees that a sequence of random variables `oscillates slowly'.

The importance of Anscombe's condition is reflected by the following result, originally obtained by Anscombe in \cite{A52}, which roughly states that if the random sequence $(N_n)_n$ is eventually close to a deterministic sequence $(k_n)_n$ which tends to $\infty$, and $(\xi_n)_n$ satisfies Anscombe's condition, then weak convergence of $\left(\xi_{n}\right)_n$ to $\xi$ implies weak convergence of $(\xi_{N_n})_n$ to $\xi$.

\begin{thm}[Anscombe]\label{thm:AnsThm}
Let $\xi$ and $(\xi_n)_n$ be $X$-valued random variables and $(N_n)_n$ $\mathbb{N}$-valued random variables. Suppose, in addition, that
\begin{itemize}
 \item[(a)] there exists $(k_n)_n$  in $\mathbb{R}^+_0$ such that $k_n \rightarrow \infty$ and $\frac{N_n}{k_n} \stackrel{\mathbb{P}}{\rightarrow} 1$,
  \item[(b)] $(\xi_n)_n$ satisfies Anscombe's condition.
 \end{itemize}
Then $\xi_n \stackrel{w}{\rightarrow} \xi \Rightarrow \xi_{N_n} \stackrel{w}{\rightarrow} \xi$.
\end{thm}

A beautiful review of the rich history of Ancombe theory is given in \cite{G12}. We wish to point out that condition (a) in Anscombe's theorem actually `kills' the randomness of $(N_n)_n$, which makes the result often inapplicable (see e.g. \cite{BLS13}, subsection 7.1.1, paragraph 2). In this paper, we will deal with this issue by using approach theory. More precisely, we will establish an approach theoretic version of Anscombe's theorem, which enhances its applicational power, in section 2. In section 3, we will apply this result to an important example in biostatistics.

\section{An approach theoretic version of Anscombe's theorem}

Define, for any sequence $(\xi_n)_n$ of $X$-valued random variables, the Anscombe index as
\begin{displaymath}
\chi_{\textrm{\upshape{Ansc}}}\left\langle(\xi_n)_n\right\rangle = \sup_{\epsilon > 0} \myinf_{\delta > 0} \limsup_{n \rightarrow \infty} \mathbb{P}\left[\max_{m = \lceil (1 - \delta) n \rceil}^{\lfloor (1 + \delta) n \rfloor} d(\xi_n, \xi_m) \geq \epsilon\right].
\end{displaymath}
Then $\chi_{\textrm{\upshape{Ansc}}}\left\langle(\xi_n)_n\right\rangle$ takes values between $0$ and $1$ and it is $0$ if and only if $(\xi_n)_n$ satisfies Anscombe's condition.

Furthermore, recall that $\lambda_w$ and $\lambda_{\mathbb{P}}$ are the limit operators for respectively the weak approach structure and the approach structure of convergence in probability (\cite{L15}). That is, for $X$-valued random variables $\xi$ and $(\xi_n)_n$,  
\begin{eqnarray}
\lefteqn{\lambda_{w}(\xi_n \rightarrow \xi)}\nonumber\\
&=& \sup_{f \in \mathcal{C}(X,\left[0,1\right])} \limsup_{n \rightarrow \infty} \left|\mathbb{E}\left[f(\xi) - f(\xi_n)\right]\right|\\
&=& \sup_{\alpha \in \mathbb{R}^+_0} \limsup_{n \rightarrow \infty} \sup_{A \textrm{ Borel }} \left(\mathbb{P}[\xi \in A] - \mathbb{P}\left[\xi_n \in A^{(\alpha)}\right]\right)\\
&=& \sup_{G \textrm{ open }} \limsup_{n \rightarrow \infty} \left(\mathbb{P}[\xi \in G] - \mathbb{P}[\xi_n \in G]\right)\\
&=& \sup_{F \textrm{ closed }} \limsup_{n \rightarrow \infty} \left(\mathbb{P}[\xi_n \in F] - \mathbb{P}[\xi \in F]\right)\\
&=& \sup_{A \in \mathcal{A}_\xi} \limsup_{n \rightarrow \infty} \left|\mathbb{P}[\xi \in A] - \mathbb{P}[\xi_n \in A]\right|,\label{eq:characLambdaW}
\end{eqnarray}
with $\mathcal{C}(X,\left[0,1\right])$ the set of continuous maps of $X$ into $[0,1]$, $A^{(\alpha)}$ the set of points from which the distance to $A$ is smaller than or equal to $\alpha$, and $\mathcal{A}_\xi$ the collection of Borel sets $A \subset X$ for which the $\mathbb{P}$-probability that $\xi$ is on the boundary of $A$ is $0$, 
and 
\begin{equation*}
\lambda_{\mathbb{P}}(\xi_n \rightarrow \xi) = \sup_{\epsilon > 0} \limsup_{n \rightarrow \infty} \mathbb{P}\left[d(\xi,\xi_n) \geq \epsilon\right].
\end{equation*}
Notice that $\lambda_{w}(\xi_n \rightarrow \xi)$ (respectively $\lambda_{\mathbb{P}}(\xi_n \rightarrow \xi)$) is a number between 0 and 1 which measures how far $(\xi_n)_n$ deviates from converging weakly (respectively in probability) to $\xi$.

In this framework, the following result holds.

\begin{thm}[Approach theoretic Anscombe]\label{thm:QAT}
Let $\xi$ and $(\xi_n)_n$ be $X$-valued random variables and $(N_n)_n$ $\mathbb{N}$-valued random variables. Then
\begin{equation}
\lambda_{w}\left(\xi_{N_n} \rightarrow \xi\right) \leq \lambda_w\left(\xi_n \rightarrow \xi\right)+ \chi_{\textrm{\upshape{Ansc}}}\left\langle(\xi_{n})_n\right\rangle + \myinf_{(k_n)_n} \lambda_{\mathbb{P}}\left(\frac{N_n}{k_n} \rightarrow 1\right),\label{QuaAns}
\end{equation}
the infimum being taken over all $(k_n)_n$ in $\mathbb{R}^+_0$ with $k_n \rightarrow \infty$.
\end{thm}

\begin{proof}
Let, for positive numbers $\zeta$, $\eta$ and $\theta$, and $(k_n)_n$ in $\mathbb{R}^+_0$ with $k_n \rightarrow \infty$,
$\lambda_w(\xi_n \rightarrow \xi) < \zeta$, 
$\chi_{\textrm{\upshape{Ansc}}}\left\langle(\xi_{n})_n\right\rangle < \eta$, and
$\lambda_{\mathbb{P}}\left(\frac{N_n}{k_n} \rightarrow 1\right) < \theta.$ Fix $\epsilon > 0$ and $F \subset X$ closed. Then there exist $\delta > 0$ and $n_1$ such that for all $n \geq n_1$ 
\begin{displaymath}
\mathbb{P}\left[\max_{m = \lceil (1 - \delta) n \rceil}^{\lfloor (1 + \delta) n \rfloor} d(\xi_n,\xi_m) \geq \epsilon\right] < \eta. 
\end{displaymath}
Also, there exists $n_2$ such that $k_n \geq n_1$ for all $n \geq n_2$.
Furthermore, there exists $n_3$ such that for all $n \geq n_3$
\begin{displaymath}
\mathbb{P}[\left|k_n - N_n\right| > \delta k_n] < \theta.
\end{displaymath}
Put $n_0 = \max\{n_1,n_2,n_3\}$. Now assume without loss of generality that $(k_n)_n$ is integer-valued. Then, for $n \geq n_0$,
\begin{eqnarray*}
\lefteqn{\mathbb{P}\left[\xi_{N_n} \in F\right] - (\eta + \theta)}\\
&\leq& \mathbb{P}\left[\left\{\xi_{N_n} \in F\right\} \cap \left\{\max_{m = \lceil (1 - \delta) k_n \rceil}^{\lfloor (1 + \delta) k_n \rfloor} d(\xi_{k_n},\xi_m) < \epsilon\right\} \cap \left\{ \left|k_n - N_n\right| \leq \delta k_n\right\} \right] \\
&\leq& \mathbb{P}\left[\xi_{k_n} \in F^{(\epsilon)}\right],
\end{eqnarray*}
whence
\begin{displaymath}
\limsup_{n \rightarrow \infty} \mathbb{P}\left[\xi_{N_n} \in F\right] \leq \limsup_{n \rightarrow \infty} \mathbb{P}\left[\xi_{k_n} \in F^{(\epsilon)}\right] + \eta + \theta \leq \mathbb{P}\left[\xi \in F^{(\epsilon)}\right] + \zeta + \eta + \theta. 
\end{displaymath}
This finishes the proof of (\ref{QuaAns}).
\end{proof}

Notice that Theorem \ref{thm:QAT} is a strict generalization of Anscombe's theorem (Theorem \ref{thm:AnsThm}). Indeed, if the conditions in Anscombe's theorem are fulfilled, then the right-hand side in (\ref{QuaAns}) becomes $0$, whence $\lambda_{w}\left(\xi_{N_n} \rightarrow \xi\right) = 0$, and we conclude that $\xi_{N_n} \stackrel{w}{\rightarrow} \xi$. Theorem \ref{thm:QAT} has the advantage that even in cases where the conditions in Anscombe's theorem fail to be fulfilled, it continues to provide valuable results. An important application in biostatistics will be treated in the next section.

\section{An application in biostatistics}

Estimation after a group sequential trial is an important theory in biostatistics. One of the key features in this theory is that the data collector is allowed to have intermediate looks at the data. After each intermediate look he can decide, based on the data observed so far and a prescribed stopping rule, whether the trial is stopped or continued. From a probabilistic point of view, randomness is imposed on the final sample size in this setting. This idea is turned into a mathematical model as follows (\cite{EF90},\cite{HP88},\cite{W92}).

Suppose that we are given a sequence $X_1, X_2, \ldots, X_n, \ldots$ of independent and identically distributed observations with mean $0$ and variance $1$. Furthermore, for each $n$, we consider a random variable $N_n$ such that 

(1) $N_n$ takes the values $n$ or $2n$, 

(2) $N_n$ is independent of $X_{n+1}, X_{n+2}, \ldots$,

(3) the conditional law 
\begin{equation}
\mathbb{P}[N_n = n \mid K_n = k] = 1_{\{\left|\cdot\right| \geq C \sqrt{n}\}}(k) \label{eq:StoppingRule}
\end{equation}
holds, where $K_n = \sum_{i = 1}^n X_i$, $C \in \mathbb{R}^+_0$ is a fixed constant, and $1_{\{\left|\cdot\right| \geq C \sqrt{n}\}}$ is the characteristic function of the set $\{x \in \mathbb{R} \mid \left|x\right| \geq C \sqrt{n}\}$.
Thus, for fixed $n$, an intermediate look is taken after having collected the data $X_1, \ldots, X_n$, and, if the sum of the data returns an extreme value in the sense that $\left|K_n\right| \geq C \sqrt{n}$, the trial is stopped, i.e. the final sample size is $N_n = n$, and, if $\left|K_n\right| < C \sqrt{n}$, the trial is continued, i.e. the additional data $X_{n+1}, \ldots, X_{2n}$ are collected, and the final sample size is $N_n = 2n$.

The majority of the existing literature on group sequential trials states that, in the above setting, conventional estimators for the mean, such as the sample average $\widehat{\mu}_{N_n} = K_{N_n}/N_n$, fail to be correct, as stopping rules generally cause them to have a large bias, a large mean squared error, and unreliable confidence intervals. Therefore, a variety of ad-hod estimation procedures has been proposed (\cite{BLS13}).

However, very recently, building upon \cite{MKA14}, a general result was obtained in \cite{BIM}, from which it can be deduced that, assuming that the $X_i$ are normally distributed in the above setting, the bias $\mathbb{E}[\widehat{\mu}_{N_n}]$ vanishes with rate $1/\sqrt{n}$, and the mean squared error $\mathbb{E}[\widehat{\mu}_{N_n}^2]$ vanishes with rate $1/n$. Also, these rates were shown to be optimal. But asymptotic normality and confidence intervals turned out to be much harder to analyze. 

Here we will show that the approach theoretic version of Anscombe's theorem (Theorem \ref{thm:QAT}) allows us to deal with asymptotic normality and confidence intervals for the estimator $\widehat{\mu}_{N_n}$, even without having to impose the additional assumption that the $X_i$ have a normal distribution.

Let $X$ be a standard normally distributed random variable.

\begin{lem}\label{lem:1}
The sequence $\left(\frac{1}{\sqrt{n}} \sum_{i=1}^n X_i\right)_n$ converges weakly to $X$. In particular, 
$$\lambda_w\left(\frac{1}{\sqrt{n}} \sum_{i=1}^n X_i \rightarrow X\right) = 0.$$ 
\end{lem}

\begin{proof}
This is a straightforward application of the central limit theorem.
\end{proof}

\begin{lem}\label{lem:2}
The sequence $\left(\frac{1}{\sqrt{n}} \sum_{i = 1}^n X_i\right)_n$ satisfies Anscombe's condition. In particular, 
$$\chi_{\textrm{\upshape{Ansc}}}\left\langle\left(\frac{1}{\sqrt{n}} \sum_{i=1}^n X_i\right)_n\right\rangle = 0.$$
\end{lem}

\begin{proof}
This follows from Kolmogorov's inequality (see e.g. the proof of Theorem 3.1 in \cite{G88}, p.16).
\end{proof}

\begin{lem}\label{lem:3}
We have
\begin{equation}
\myinf_{(k_n)_n} \lambda_{\mathbb{P}}\left(\frac{N_n}{k_n} \rightarrow 1\right) = \min \{\mathbb{P}\left[- C < X < C\right],1 - \mathbb{P}\left[- C < X < C\right] \},\label{eq:noAns}
\end{equation}
the infimum being taken over all $(k_n)_n$ in $\mathbb{R}^+_0$ with $k_n \rightarrow \infty$, and $C$ being the fixed constant determining the stopping rule (\ref{eq:StoppingRule}).
\end{lem}

\begin{proof}
The conditional law 
$$\mathbb{P}[N_n = n \mid K_n = k] = 1_{\{\left|\cdot\right| \geq C \sqrt{n}\}}(k)$$
implies that
\begin{equation}
\mathbb{P}[N_n = n] = \mathbb{P}\left[\left|K_n\right| \geq C \sqrt{n}\right] = 1 - \mathbb{P}\left[- C < K_n/\sqrt{n} < C\right].\label{eq:ProbNisn}
\end{equation}
Thus
\begin{equation}
\mathbb{P}[N_n = 2n] = 1 - \mathbb{P}[N_n = n] = \mathbb{P}\left[- C < K_n/\sqrt{n} < C\right].\label{eq:ProbNis2n}
\end{equation}
Furthermore, taking $\epsilon = 1/3$, we see that, for any sequence $(k_n)_n$ in $\mathbb{R}^+_0$, 
$$\left|\frac{n}{k_n} - 1\right| < \epsilon \Leftrightarrow \frac{3n}{4} < k_n < \frac{3n}{2}$$
and
$$\left|\frac{2n}{k_n} - 1\right| < \epsilon \Leftrightarrow \frac{3n}{2} < k_n < \frac{n}{3},$$
from which it follows that at least one of the events
\begin{equation*}
A_n = \left\{\left|\frac{N_n}{k_n} - 1\right| < \epsilon, N_n = n\right\}
\end{equation*}
and
\begin{equation*}
A_{2n} = \left\{\left|\frac{N_n}{k_n} - 1\right| < \epsilon, N_n = 2n\right\}
\end{equation*}
must be empty. Observe that
\begin{eqnarray*}
\lefteqn{\mathbb{P}\left[\left|\frac{N_n}{k_n} - 1\right| \geq \epsilon\right]}\\
&=& 1 -  \mathbb{P}\left[\left|\frac{N_n}{k_n} - 1\right| < \epsilon\right]\\
&=& 1 - \mathbb{P}\left[\left|\frac{N_n}{k_n} - 1\right| < \epsilon, N_n = n\right] - \mathbb{P}\left[\left|\frac{N_n}{k_n} - 1\right| < \epsilon, N_n = 2n\right].
\end{eqnarray*}
Thus, if $A_n = \emptyset$, by (\ref{eq:ProbNis2n}),
\begin{eqnarray*}
\mathbb{P}\left[\left|\frac{N_n}{k_n} - 1\right| \geq \epsilon\right] &=& 1 - \mathbb{P}\left[\left|\frac{N_n}{k_n} - 1\right| < \epsilon, N_n = 2n\right]\\
 &\geq& 1 - \mathbb{P}[N_n = 2n]\\
 &=&  \mathbb{P}\left[- C < K_n/\sqrt{n} < C\right],
\end{eqnarray*}
and, if $A_{2n} = \emptyset$, by (\ref{eq:ProbNisn}),
\begin{eqnarray*}
\mathbb{P}\left[\left|\frac{N_n}{k_n} - 1\right| \geq \epsilon\right] &=& 1 - \mathbb{P}\left[\left|\frac{N_n}{k_n} - 1\right| < \epsilon, N_n = n\right]\\
 &\geq& 1 - \mathbb{P}[N_n = n]\\
 &=&   1 - \mathbb{P}\left[- C < K_n/\sqrt{n} < C\right].
\end{eqnarray*}
Therefore, 
\begin{eqnarray*}
\lefteqn{\mathbb{P}\left[\left|\frac{N_n}{k_n} - 1\right| \geq \epsilon\right]}\\
&\geq& \min \left\{\mathbb{P}\left[- C < K_n/\sqrt{n} < C\right],1 - \mathbb{P}\left[- C < K_n/\sqrt{n} < C\right]\right\},
\end{eqnarray*}
whence, using the fact that $K_n/\sqrt{n} \stackrel{w}{\rightarrow} X$ by the central limit theorem, 
\begin{eqnarray}
\lefteqn{\myinf_{(k_n)_n} \sup_{\epsilon > 0} \limsup_{n \rightarrow \infty} \mathbb{P}\left[\left|\frac{N_n}{k_n} - 1\right| \geq \epsilon\right]} \label{eq:noAnslow}\\
&\geq& \min \{\mathbb{P}\left[- C < X < C\right],1 - \mathbb{P}\left[- C < X < C\right] \}.\nonumber
\end{eqnarray}
On the other hand, choosing $k_n = n$, gives, for $\epsilon$ small, by (\ref{eq:ProbNisn}),
\begin{equation*}
\mathbb{P}\left[\left|\frac{N_n}{k_n} - 1\right| \geq \epsilon\right] = \mathbb{P}\left[- C < K_n/\sqrt{n} < C\right],
\end{equation*}
and, choosing $k_n = 2n$, gives, for $\epsilon$ small, by (\ref{eq:ProbNis2n}),
\begin{equation*}
\mathbb{P}\left[\left|\frac{N_n}{k_n} - 1\right| \geq \epsilon\right] = 1 - \mathbb{P}\left[- C < K_n/\sqrt{n} < C\right].
\end{equation*}
That is, again by the central limit theorem,
\begin{eqnarray}
\lefteqn{\myinf_{(k_n)_n} \sup_{\epsilon > 0} \limsup_{n \rightarrow \infty} \mathbb{P}\left[\left|\frac{N_n}{k_n} - 1 \right|\geq \epsilon\right] } \label{eq:noAnsup}\\
&\leq& \min \{\mathbb{P}\left[- C < X < C\right],1 - \mathbb{P}\left[- C < X < C\right] \}.\nonumber
\end{eqnarray}
Combining (\ref{eq:noAnslow}) and (\ref{eq:noAnsup}), proves (\ref{eq:noAns}).
\end{proof}

Notice that if condition (a) in Anscombe's theorem (Theorem \ref{thm:AnsThm}) is satisfied, then expression (\ref{eq:noAns}) vanishes, which never happens. Therefore, the previous lemma shows that Anscombe's theorem fails to be applicable in this setting. However, without additional effort, we can continue to use the more powerful approach theoretic version of Anscombe's theorem (Theorem \ref{thm:QAT}) to obtain the following result.

\begin{thm}
The inequality
\begin{equation}
\lambda_{w}\left(K_{N_n}/\sqrt{N_n} \rightarrow X\right) \leq \min \{\mathbb{P}\left[- C < X < C\right],1 - \mathbb{P}\left[- C < X < C\right] \}\label{eq:MainIneq}
\end{equation}
holds, with $C$ the fixed constant determining the stopping rule (\ref{eq:StoppingRule}).
\end{thm}

\begin{proof}
Apply Theorem \ref{thm:QAT} and notice that in the upper bound in (\ref{QuaAns}) the first term vanishes by Lemma \ref{lem:1}, the second term vanishes by Lemma \ref{lem:2}, and the third term equals $\min \{\mathbb{P}\left[- C < X < C\right],1 - \mathbb{P}\left[- C < X < C\right] \}$ by Lemma \ref{lem:3}.
\end{proof}

Using (\ref{eq:characLambdaW}), we can conclude from (\ref{eq:MainIneq}) that for each Borel set $A \subset \mathbb{R}$ for which $X$ has zero probability of being on the boundary of $A$, and each $\epsilon > 0$, there exists $n_0$ such that 
\begin{eqnarray*}
\lefteqn{\left|\mathbb{P}\left[K_{N_n}/\sqrt{N_n} \in A\right] - \mathbb{P}[X \in A] \right| }\\
&\leq& \min \{\mathbb{P}\left[- C < X < C\right],1 - \mathbb{P}\left[- C < X < C\right]\} + \epsilon
\end{eqnarray*}
for all $n \geq n_0$. In particular, recalling that $\widehat{\mu}_{N_n} = K_{N_n}/N_n$, and taking $A = [-B,B]$ with $B \in \mathbb{R}^+_0$, we infer that for each $\epsilon > 0$ there exists $n_0$ such that
\begin{eqnarray*}
\lefteqn{\left|\mathbb{P}\left[\widehat{\mu}_{N_n} -B/\sqrt{N_n} \leq 0 \leq \widehat{\mu}_{N_n} + B/\sqrt{N_n}\right] - \mathbb{P}[- B \leq X \leq B] \right| }\\
&\leq& \min \{\mathbb{P}\left[- C < X < C\right],1 - \mathbb{P}\left[- C < X < C\right]\} + \epsilon
\end{eqnarray*}
for all $n \geq n_0$. We infer that the approach theoretic version of Anscombe's theorem justifies the use of confidence intervals based on the estimator $\widehat{\mu}_{N_n}$ if $n$ is sufficiently large and $\min\{\mathbb{P}\left[- C < X < C\right],1 - \mathbb{P}\left[- C < X < C\right]\}$ is sufficiently small.


\begin{thebibliography}{99}
\bibitem[A52]{A52} Anscombe, F. J. (1952) {\em Large-sample theory of sequential estimation.} Proc. Cambridge Philos. Soc. 48, 600--607.
\bibitem[BLS13]{BLS13} Bartroff, J.; Lai, T. L.; Shih, M.-C. (2013) {\em Sequential experimentation in clinical trials. Design and analysis.} Springer Series in Statistics. Springer, New York.
\bibitem[BIM]{BIM} Berckmoes, B.; Ivanova, A.; Molenberghs, G. (2017) {\em On the sample mean after a group sequential trial} https://arxiv.org/abs/1706.01291.
\bibitem[EF90]{EF90} 
Emerson, S. S.; Fleming, T. R. (1990).
{\em Parameter estimation following group sequential hypothesis testing.}
Biometrika 77, 875--892.
\bibitem[G88]{G88} Gut, A. (1988) {\em Stopped random walks. Limit theorems and applications.} Applied Probability. A Series of the Applied Probability Trust, 5. Springer-Verlag, New York.
\bibitem[G12]{G12} Gut, A. (2012) {\em Anscombe's theorem 60 years later.} Sequential Anal. 31, no. 3, 368--396.
\bibitem[HP88]{HP88}
Hughes, M.D.; Pocock, S.J. (1988).
{\em Stopping rules and estimation problems in clinical trials.}
Statistics in Medicine 7, 1231--1242.
\bibitem[L15]{L15} Lowen, R. (2015) {\em Index analysis. Approach theory at work.} Springer Monographs in Mathematics. Springer, London.
 \bibitem[MKA14]{MKA14} Molenberghs, G.; Kenward, M. G.; Aerts, M.; Verbeke, G.; Tsiatis, A. A.; Davidian, M.; Rizopoulos, D. (2014) {\em On random sample size, ignorability, ancillarity, completeness, separability, and degeneracy: sequential trials, random sample sizes, and missing data.} Stat. Methods Med. Res. 23, no. 1, 11--41. 
\bibitem[W92]{W92}
Woodroofe, M. (1992)
{\em Estimation after sequential testing: a simple approach for a truncated sequential probability ratio test.} 
Biometrika 79, no. 2, 347--353.
\end{thebibliography}
\end{document}